\documentclass[preprint,12pt]{elsarticle}

\usepackage{amsmath}
\usepackage{amssymb}
\usepackage{graphicx}
\usepackage{color} 
\usepackage[normalem]{ulem}
\usepackage[utf8]{inputenc}

\usepackage[
centering,
marginparwidth=0in,
textwidth=6.5in,
marginparsep=2em,
top=2.5cm,
bottom=2.5cm,
]{geometry}

\usepackage{graphicx}% Include figure files
\usepackage{dcolumn}% Align table columns on decimal point
\usepackage{bm}% bold math
\usepackage{amsmath}
\usepackage{amssymb}
\usepackage{mathtools}
\usepackage{epstopdf, epsfig}
\usepackage{color}
\usepackage{tikz}
\usepackage{pgfplots}

\newcommand{\surf}{\mathcal{S}}

\newcommand{\R}{\mathbb{R}}

\newcommand{\normal}{\boldsymbol{\nu}}
\newcommand{\normalC}{\nu}

\newcommand{\vel}{\boldsymbol{u}}
\newcommand{\velC}{u}
\newcommand{\velTest}{\boldsymbol{\phi}}
\newcommand{\relvel}{\boldsymbol{v}}

\newcommand{\param}{\boldsymbol{X}}
\newcommand{\paramC}{X}
\newcommand{\paramTest}{\boldsymbol{Z}}
\newcommand{\paramTestC}{Z}
\newcommand{\paramUpdate}{\boldsymbol{Y}}
\newcommand{\paramUpdateC}{{Y}}
\newcommand{\coord}{\boldsymbol{x}}
\newcommand{\coordC}{x}

\newcommand{\shop}{\boldsymbol{B}}
\newcommand{\shopC}{B}
\newcommand{\meanc}{\mathcal{H}}

\newcommand{\proj}{\boldsymbol{P}}
\newcommand{\Id}{\boldsymbol{I}}

\newcommand{\DivP}{\operatorname{div}_{\!\proj}\!}

\newcommand{\inputTikzPic}[1]{\ifthenelse{\boolean{plotTikzPics}}{\input{#1}}{\fbox{\centering\begin{minipage}[t][5cm][t]{0.9\textwidth}\centering\textbf{\color{red}enable plotTikzPics}\end{minipage}}}}

\newcommand{\Grad}{\nabla}
\newcommand{\Div}{\operatorname{div}}%
\newcommand{\GradSurf}{\Grad_{\!\surf}}
\newcommand{\GradP}{\Grad_{\!\proj}}

\newcommand{\Be}{\textup{Be}}

\newcommand{\ReynoldsNumber}{\mathrm{Re}}
\newcommand{\BendingNumber}{\mathrm{Be}}

\renewcommand{\Re}{\ReynoldsNumber}
\renewcommand{\Be}{\BendingNumber}
\newcommand{\stressC}{\sigma}
\newcommand{\stress}{\boldsymbol{\stressC}}
\newcommand{\float}[1]{\mathrm{e}\!-\!{#1}}

\begin{document}

\begin{frontmatter}
\title{A numerical approach for fluid deformable surfaces with conserved enclosed volume}

\author[label1]{Veit Krause}
\affiliation[label1]{organization={Institut f\"ur wissenschaftliches Rechnen, TU Dresden},
            city={Dresden},
            postcode={01062}, 
            country={Germany}}

\author[label1,label2,label3]{Axel Voigt}
\affiliation[label2]{organization={Center for Systems Biology Dresden (CSBD)},
            addressline={Pfotenhauerstra{\ss}e 108}, 
            city={Dresden},
            postcode={01307}, 
            country={Germany}}
\affiliation[label3]{organization={Cluster of Excellence Physics of Life (PoL)},
            city={Dresden},
            postcode={01062}, 
            country={Germany}}

\begin{abstract}
We consider surface finite elements and a semi-implicit time stepping scheme to simulate fluid deformable surfaces. Such surfaces are modeled by incompressible surface Navier-Stokes equations with bending forces. Here, we consider closed surfaces and enforce conservation of the enclosed volume. The numerical approach builds on higher order surface parameterizations, a Taylor-Hood element for the surface Navier-Stokes part, appropriate approximations of the geometric quantities of the surface mesh redistribution and a Lagrange multiplier for the constraint. The considered computational examples highlight the solid-fluid duality of fluid deformable surfaces and demonstrate convergence properties that are known to be optimal for different sub-problems.
\end{abstract}

\end{frontmatter}

\section{Introduction}

Fluid deformable surfaces are thin fluidic sheets of soft materials exhibiting a solid–fluid duality.
While they store elastic energy when stretched or bent, as solid shells, under in-plane shear, they flow as viscous two-dimensional fluids. This duality has several consequences: it establishes a tight interplay between tangential flow and surface deformation. In the presence of curvature, any shape change is accompanied by tangential flow and, vice versa, the surface deforms due to tangential flow \cite{torres2019modelling,Voigt_JFM_2019,reuther2020numerical}. Surfaces with these properties are ubiquitous interfaces in biology, playing essential roles in processes from the subcellular to the tissue scale, e.g., describing the dynamics of cortical flow during cell division \cite{Mayeretal_Nature_2010} and the shape evolution of epithelia monolayers in morphogenesis \cite{Heisenbergetal_Cell_2013}. 

The mathematical description of fluid deformable surfaces has been introduced in \cite{Arroyoetal_PRE_2009,Salbreuxetal_PRE_2017,torres2019modelling,reuther2020numerical}. The model also relates to two-phase flow problems with Boussinesq–Scriven interface condition \cite{Nitschkeetal_JFM_2012,Barrettetal_PRE_2015} in the limit of a large Saffman–Delbrück number \cite{Saffmanetal_PNAS_1975}, where bulk and surface flow can be decoupled. Numerical realizations of these models, which are surface Navier-Stokes equations on evolving surfaces, are restricted to special cases. In \cite{Arroyoetal_PRE_2009,Mietkeetal_PNAS_2019,torres2019modelling,Sahuetal_JCP_2020,reuther2020numerical} an axisymmetric setting is used, the approach is restricted to simply-connected surfaces, a Stokes approximation is considered or constraints on the evolution are omitted. Furthermore, non of the previous numerical attempts considered higher order methods. In view of the recently discovered numerical analysis results for a surface vector-Laplacian \cite{Grossetal_SIAMJNA_2018,Hansboetal_IMAJNA_2020,hardering2021tangential} and a surface Stokes problem \cite{Olshanskiietal_SIAMJSC_2018,hardering2022Stokes} on stationary surfaces and various corresponding computational studies \cite{Fries_IJNMF_2018,Reutheretal_PF_2018,Grossetal_JCP_2018,nestler2019finite,Ledereretal_IJNME_2020,Ranketal_PF_2021,Brandneretal_SIAMJSC_2022,Bachinietal_arXiv_2022}, which require higher order approximations of geometric properties of the surface to guarantee convergence, previous numerical approaches \cite{torres2019modelling,reuther2020numerical} might not converge to the desired solution. We here close this gap and provide a numerically feasible approach for the full problem which experimentally shows higher order (presumably optimal) convergence properties in space and is shown to be applicable in a general context.

The paper is structured as follows: In Section \ref{sec:Context} we introduce the surface differential operators and the continuous model based on \cite{jankuhn2018,torres2019modelling,reuther2020numerical}. Moreover we combine the surface Navier-Stokes model with bending forces and conservation of enclosed volume with an approach of \cite{barrett2008parametric} to locally redistribute mesh points which improves the mesh quality during evolution without influencing the dynamics. In Section \ref{sec:numerics} we introduce a higher order polynomial approximation of the discrete surface \cite{praetorius2020dunecurvedgrid} and a surface finite element method (SFEM) for vector-valued functions \cite{nestler2019finite}. The discretization builds on the SFEM approach for the surface Stokes equation on stationary surfaces used in \cite{Brandneretal_SIAMJSC_2022}, the discretization for Willmore flow considered in \cite{Dziuk_NM_2008} and a volume constraint considered by a Lagrange multiplier. Numerical results are shown in Section \ref{sec:results} and are chosen to highlight the solid-fluid duality and the influence of critical parameters, such as the Reynolds number $\Re$ and the bending capillary number $\Be$. Due to a lack of numerical analysis results for the full problem we consider experimental convergence studies. For the fluid properties we demonstrate that optimal order of convergence known for the Stokes problem on stationary surfaces can also been achieved for the full problem. Also for the solid properties we achieve the same order of convergence for the full problem, which had been obtained for pure Willmore flow without any constraints. The increased accuracy leads to qualitative differences compared with \cite{reuther2020numerical} if the volume constraint is dropped. We validate the obtained solution with semi-analytic results assuming an axisymmetric setting. We further demonstrate the applicability of the mesh redistribution algorithm. For extreme shape changes, as in the deforming Stanford bunny, the considered local approach reaches its limits but still provides sufficiently regular meshes to solve the linear systems. Finally we summarize our results and draw conclusions in Section \ref{sec:Discussion}. 

\section{Continuous Problem}
\label{sec:Context}

We consider a time dependent smooth domain $\Omega=\Omega(t)\subset\R^3$ and its boundary $\surf = \surf(t)$, which is given by a parametrization $\param$. Related to $\surf$ we denote the surface normal $\normal$, the shape Operator $\shop$ with $\shopC_{ij} = -\GradSurf^i \normalC_j$, the mean curvature $\meanc= \operatorname{tr}\shop$ and the surface projection $\proj=\Id-\normal\otimes\normal$. Let $\GradSurf$ and $\Div_{\surf}$ be the surface differential operators with respect to the covariante derivative. Those operators are well defined for vector fields in the tangent bundle of $\surf$. For a non-tangential vector field $\vel\colon \surf\rightarrow\R^3$ we use the tangential derivative $\GradP\vel = \proj\nabla \vel^e\proj$ and the tangential divergence $\DivP \vel = \operatorname{tr}[\proj\nabla \vel^e]$ where $\vel^e$ is an extension of $\vel$ constant in normal direction and $\nabla$ the gradient of the embedding space $\R^3$. $\GradP\vel$ is a pure tangential tensor field and relates to the covariante operator by $\GradP\vel = \GradSurf(\proj\vel)-(\vel\cdot\normal)\shop$. Similarly, it hold $\DivP\vel = \Div_{\surf}(\proj\vel) - (\vel\cdot\normal)\meanc$, see \cite{Brandneretal_SIAMJSC_2022}.

We consider the incompressible surface Navier-Stokes model on $\surf$, independently derived in \cite{Arroyoetal_PRE_2009,Reutheretal_MMS_2015,Kobaetal_QAM_2017,Miura_QAM_2018,jankuhn2018,Nitschkeetal_PRF_2019,torres2019modelling,reuther2020numerical,Sahuetal_JCP_2020} and previously considered numerically in \cite{torres2019modelling,reuther2020numerical}. The derivations differ with respect to the physical principles used and in the coordinate system in which the resulting equations are represented. The different formulations have been compared in \cite{Brandneretal_arXiv_2021} and at least with the corrections \cite{Reutheretal_MMS_2018,Kobaetal_QAM_2018} and \cite{Yavarietal_JNS_2016} concerning the model in \cite{Arroyoetal_PRE_2009}, all agree and can be formulated as:
\begin{align}
    \label{eq:navierstokes}
    \begin{aligned}
        \partial_t \vel + \nabla_{\relvel}\vel &= -\GradSurf p-p\meanc\normal + \frac{2}{\Re}\DivP \stress + \mathbf{b}  \\
        \DivP \vel &= 0  
    \end{aligned}    
\end{align}
where $\vel$ is the surface fluid velocity, $p$ the surface pressure, $\stress = \frac{1}{2}(\GradP\vel + \GradP\vel^T)$ the rate of deformation tensor, $\Re$ the Reynolds number and $\mathbf{b}$ an acting force. The convection term is defined by $[\nabla_{\relvel}\vel]_i=\GradSurf\velC_i\cdot\relvel$, $i=0,1,2$ where $\relvel$ is the relative material velocity of the surface. The surface fluid velocity is related to the parametrization by $(\vel\cdot\normal)\normal = \partial_t \param$. This formulation thus provides an Eulerian approach in tangential direction but a Lagrangian approach in normal direction.

Eqs. \eqref{eq:navierstokes} imply conservation of the surface area $\vert\surf\vert$, as a consequence of the local inextensibility constraint. We in addition want to enforce conservation of the enclosed volume $\vert\Omega\vert$. This can be realized by $\int_\surf \vel\cdot\normal = 0$ and incorporated in eqs. \eqref{eq:navierstokes} by a scalar Lagrange multiplier $\lambda$, as explained in detail in \cite{torres2019modelling}. We get  
\begin{align}
    \label{eq:conservNavierStokes}
    \begin{aligned}
        \partial_t \vel + \nabla_{\relvel}\vel &= -\GradSurf p-p\meanc\normal + \frac{2}{\Re}\DivP \stress-\lambda\normal + \mathbf{b} \\
        \DivP \vel &= 0\\
        \int_\surf \vel\cdot\normal \; \mathrm{d}\surf &= 0.
    \end{aligned}
\end{align}
As in \cite{torres2019modelling,reuther2020numerical} we assume elastic properties of the surface $\surf$ if deformed in normal direction and define $\mathbf{b}$ as the contribution resulting from the Helfrich energy \cite{Helfrich+1973+693+703}, $E_{H}=\frac{1}{2\Be}\int_\surf (\meanc-\meanc_0)^2 \mathrm{d}\surf$, with the bending capillary number $\Be$ and the spontaneous curvature $\meanc_0$. We will only consider the case $\meanc_0=0$ and obtain
\begin{align}
    \mathbf{b} &= \frac{1}{\Be}(-\Delta_\surf\meanc -\meanc(\Vert\shop\Vert^2 - \frac{1}{2}\operatorname{tr}(\shop)^2)\normal. \label{eq:bending}
\end{align}

Using this Eulerian-Lagrangian description for eqs. \eqref{eq:conservNavierStokes} and \eqref{eq:bending} also in the discrete formulation, to be considered below, the regularity of the underlying mesh cannot be guaranteed. To cure this issue we follow an approach introduced in \cite{barrett2008parametric} by adding a tangential flow of the mesh points. Instead of $(\vel\cdot\normal)\normal = \partial_t \param$ we consider 
\begin{align}
    \label{eq:meshflow}
    \begin{aligned}
        \partial_t \param \cdot \normal &= \vel\cdot \normal \\
        \meanc \normalC_i &= \Delta_{\surf} \paramC_i  \quad i=0,1,2.
    \end{aligned}
\end{align}
This does not alter the movement in normal direction but leads to a local redistribution of the mesh points in the discrete system in tangential direction. Moreover by eq. \eqref{eq:meshflow} we get a characterization of $\meanc$, which has been used numerically in various context \cite{Dziuk_NM_1990,Baenschetal_JCP_2005,Hausseretal_JSC_2007,Dziuk_NM_2008}. 
We thus obtain the relative material velocity $\relvel=\vel-\partial_t\param$ which enters the convective term $\nabla_{\relvel}\vel$ in eqs. \eqref{eq:navierstokes} or \eqref{eq:conservNavierStokes}. Combining eqs. \eqref{eq:conservNavierStokes} - \eqref{eq:meshflow} leads to a numerically feasible incompressible surface Navier-Stokes model for fluid deformable surfaces with conservation of the enclosed volume. 

\section{Numerical approach}
\label{sec:numerics}

\subsection{Surface approximation} 
Let $\surf_h$ be a discrete higher order approximation of $\surf$. The construction of such surfaces is explained in \cite{demlow2006localized,Brandneretal_SIAMJSC_2022} and implemented for the \textsc{Dune} finite element framework \cite{bastian2021dune}, see \cite{praetorius2020dunecurvedgrid,dunealugrid:16}. Let $\surf_h^{\mathrm{lin}} = \bigcup_{\hat{T}\in \mathcal{T}^{\mathrm{lin}}} \hat{T}$ be a piecewise linear reference surface given by a shape regular surface triangulation $\mathcal{T}^{\mathrm{lin}}$ and $\param\colon\surf_h^{\mathrm{lin}}\rightarrow\surf$ be a bijective map such that $\surf=\bigcup_{\hat{T}\in \mathcal{T}^{\mathrm{lin}}} \param(\hat{T})$. Let $I_k$ be the $k$th order interpolation such that $I_k\param = \param_k \in \mathbb{P}_k(\hat{T})^3$ for each triangle. This allows us to describe the curved mesh by $\mathcal{T}=\{T=\param_k(\hat{T}) \vert \hat{T}\in\mathcal{T}^{\mathrm{lin}}\}$ and the discrete surface is given by 
\begin{align*}
    \surf_h = \bigcup_{\hat{T}\in \mathcal{T}^{\mathrm{lin}}} \param_k(\hat{T}) = \bigcup_{T\in \mathcal{T}} T.
\end{align*}
We define the mesh size as $h=h_{max}=\max_{T\in\mathcal{T}} h_T$ where the size of a triangle $h_T$ is given by the longest edge of $T$. For better readability we will use $\param$ instead of $\param_k$ and also denote the discrete normal vector of $\surf_h$ by $\normal$ and the discrete shape operator of $\surf_h$ by $\shop$. The discrete projection and surface operators are defined with respect to $\normal$ and $\shop$. Please consider that $\normal$ is multi-valued at the intersection of the mesh elements. It is therefore not recommended to use $\normal$ to compute the curvature information required for the bending force in eq. \eqref{eq:bending}. Different approaches have been proposed to deal with this issue, see e.g. \cite{Dziuk_NM_2008,Nitschkeetal_JFM_2012,Kovacsetal_NM_2019}. Here, we will obtain the mean curvature of $\surf_h$, again denoted by $\meanc$, from eq. \eqref{eq:meshflow}. 

\subsection{Discretization by surface finite elements}

We consider a surface finite element (SFEM) approach to discretize the problem in space. Such approaches are well established for scalar-valued surface partial differential equations \cite{dziuk2013finite} and have been extended to vector- and tensor-valued surface partial differential equations in \cite{nestler2019finite}. 

We define the $k$th order finite element space on $\surf_h$ by $V_k(\surf_h)=\{ v \in C^0(\surf_h) \vert \quad v\vert_{T} \in \mathbb{P}_k(T) \}$. Let $(\boldsymbol{a},\boldsymbol{b}) = \int_\surf \boldsymbol{a}\cdot \boldsymbol{b} \, \mathrm{d}\surf$ be the $L^2$ inner product on $\surf_h$ \cite{nestler2019finite,Brandneretal_SIAMJSC_2022}. We use $V_h = V_k(\surf_h)^3 \times V_{k-1}(\surf_h) \times V_{k}(\surf_h) \times V_{k}(\surf_h)^3$ as test and trial space which leads to $k$th order Taylor-Hood elements for the velocity (order $k$) and the pressure (order $k-1$) combined with $k$th order spaces for the mean curvature and bijective map. We further split the problem and consider first the Lagrange multiplier $\lambda$ as given. 

Again omitting any indices to distinguish between discrete and continuous quantities we obtain the space discretized system for a tuple $(\vel,p,\meanc,\param) \in V_h$ by 
\begin{align}
    \label{eq:conservNavierStokesSpaceDis}
    \begin{aligned}
        (\partial_t \vel + \nabla_{\relvel} \vel ,\velTest)  &= (p, \DivP\velTest) - \frac{2}{\Re} (\stress, \GradP\velTest) - ( \lambda \normal, \velTest)  \\
        & \quad +\frac{1}{\Be}[ ( \GradSurf\meanc, \GradSurf(\normal\cdot\velTest)) - (\beta\meanc\normal,\velTest)]       \\
        (\DivP \vel,q) &= 0 \\
        (\partial_t \param\cdot\normal , h) &= (\vel\cdot \normal, h) \\
        (\meanc \normalC_{i} , \paramTestC_i)  &= ( \GradSurf \paramC_i , \GradSurf \paramTestC_i )  \quad i=0,1,2
    \end{aligned}
\end{align}
for all $(\velTest,q,h,\paramTest)\in V_h$ where $\beta =\Vert\shop\Vert^2 - \frac{1}{2}\operatorname{tr}(\shop)^2$ and $\lambda$ chosen such that $\vel$ fulfills the volume constraint. This in principle requires to solve a non-local problem. Moreover we used the identity $(-\GradSurf p -p\meanc\normal,\velTest)=(p, \DivP\velTest)$. 

\subsection{Discretization in time}\label{sec:alternativeApproach}

The volume constraint makes the system of surface partial differential equations non-local, which further complicates approaches like mesh parallelization and multiprocessor computation. The following discretization considers an operator splitting approach which decouples the system of surface partial differential equation and the non-local term. It treats the Lagrange multiplier $\lambda$ as a given quantity for solving the surface partial differential equations, as in eqs. \eqref{eq:conservNavierStokesSpaceDis}, and considers a Newton iteration to update $\lambda$ to fulfill the volume constraint in each timestep. This approach allows for a straight forward parallelization.

Let $\{t^n\}_{n=0}^N$ be a set of times for timestep width $\tau>0$ with $t^n=n\tau$. Let $\surf_h^n=\surf_h(t^n)$ and each variable with subscript $n$ corresponds to the timestep $t^n$. We discretize the system in time by a semi-implicite Euler scheme. That means that all linear terms will be discretized by the implicite Euler method and all nonlinear terms are treated semi-implicitly. We use scalar products, the normal vector $\normal=\normal(t^n)$ and the shape operator $\shop=\shop(t^n)$ with respect to $\surf_h^n$. This also linearizes the bending force and all terms containing the normal vector. Instead of the surface parametrization we compute a surface update $\paramUpdate^{n+1} = \param^{n+1}-\param^n$ as in \cite{barrett2008parametric} and take the relative material velocity $\relvel^n = \vel^n - \frac{1}{\tau}\paramUpdate^n$ which also linearizes the convective term in the system. \\
We define the bilinear form $L^n\colon V_h(\surf_h^n)\times V_h(\surf_h^n)\rightarrow\R$ and the linear forms $r^n\colon V_h(\surf_h^n)\rightarrow\R$ by 
\begin{align*}
    L^n((\vel,p,,\meanc,\paramUpdate),(\velTest,q,h,\paramTest)) &= 
    \left(
    \begin{aligned}
    &(\vel + \tau \nabla_{\relvel^n} \vel ,\velTest) - (\tau p,\DivP\velTest) + \frac{2\tau}{\Re} (\stress, \GradP\velTest)\\
    &\qquad -\frac{\tau}{\Be}[( \GradSurf\meanc, \GradSurf(\normal\cdot\velTest)) - (\beta\meanc\normal,\velTest)] \\
    &(\DivP \vel,q) \\
    &(\paramUpdate\cdot\normal , h) - (\tau\vel\cdot\normal , h) \\    
    &(\meanc\normalC_i , \paramTestC) + (\GradSurf\paramUpdateC_i , \GradSurf\paramTestC_i ), \quad i=0,1,2
    \end{aligned}
    \right) \\\\
    r^n((\velTest,q,h,\paramTest)) &= 
    \left(
    \begin{aligned}
    &(\vel^n,\velTest)  \\
    & 0 \\
    & 0 \\
    &(\GradSurf \paramC^n_i , \GradSurf\paramTestC_i ), \quad i=0,1,2
    \end{aligned}
    \right)
\end{align*} 
where the index $n$ highlights the dependency of $L^n$ and $r^n$ on the time steps. We thus seek for solutions   $(\hat{\vel}^{n+1},\hat{p}^{n+1},\hat{\meanc}^{n+1},\hat{\paramUpdate}^{n+1})\in V_h(\surf_h^n)$ of the variational problem 
\begin{align} \label{eq:separate1}
    L^n((\hat{\vel}^{n+1},\hat{p}^{n+1},\hat{\meanc}^{n+1},\hat{\paramUpdate}^{n+1}),(\velTest,q,h,\paramTest)) = r^n((\velTest,q,h,\paramTest)) - ((\tau\lambda\normal,\velTest),0,0,0) 
\end{align}
$\forall (\velTest,q,h,\paramTest)\in V_h(\surf_h^n)$ with $\lambda$ assumed to be given. Thereby the symbol $\hat{\cdot}$ denotes the new solution at time $t^{n+1}$ but with respect to the old geometry $\surf_h^n$. After lifting the solution to the new geometry $\surf_h^{n+1}$ the symbol $\hat{\cdot}$ will be dropped.

In order to fulfill the volume constraint $\lambda$ has to be chosen such that 
\begin{align}\label{eq:separate2}
  \Phi(\lambda)\coloneqq\int_{\surf_h^n} \hat{\vel}^{n+1}(\lambda)\cdot\normal \mathrm{d}\surf = 0.
\end{align}
We consider $\Phi(\lambda) = 0$ as a non-linear equation in $\lambda$, which motivates to apply a Newton method.
More precisely, let $\lambda^j\in\R$ and $(\hat{\vel}^j,\hat{p}^j,\hat{\meanc}^j,\hat{\paramUpdate}^j)\in V_h(\surf_h^n)$ for $j\in\mathbb{N}$ generated by the Newton iterations, for each Newton step we compute the Newton update $\delta\lambda^{j+1}$ by $\left(\frac{d}{d\lambda}\Phi(\lambda^j)\right) \delta\lambda^{j+1} = -\Phi(\lambda^j)$ where $\Phi(\lambda^j)=\int_{\surf_h^n} \hat{\vel}^{j+1}\cdot\normal$ and $\frac{d}{d\lambda}\Phi(\lambda^j) = \int_{\surf_h^n} (\frac{d}{d\lambda}\vert_{\lambda^j}\hat{\vel}^{j+1})\cdot\normal \mathrm{d}\surf$. To compute the derivative of $\Phi$ we will insert the derivative of the solution with respect to $\lambda$ in the bilinear form
\begin{align*}
    L^n&\left(\left(
    \frac{d}{d\lambda}\Big\vert_{\lambda^j}\hat{\vel}^{j+1},
    \frac{d}{d\lambda}\Big\vert_{\lambda^j}\hat{p}^{j+1}, 
    \frac{d}{d\lambda}\Big\vert_{\lambda^j}{\meanc}^{j+1},
    \frac{d}{d\lambda}\Big\vert_{\lambda^j}\hat{\paramUpdate}^{j+1}
    \right),
    (\velTest,q,h,\paramTest)\right) \\
    &= \frac{d}{d\lambda}\Big\vert_{\lambda^j} L^n ((\hat{\vel}^{j+1},\hat{p}^{j+1} ,{\meanc}^{j+1},\hat{\paramUpdate}^{j+1}),(\velTest,q,h,\paramTest)) \\
    &= \frac{d}{d\lambda}\Big\vert_{\lambda^j}( r^n((\velTest,q,h,\paramTest)) + (-\tau \lambda^j(\normal , \velTest),0,0,0)^T)) \\
    &= (-\tau (\normal , \velTest),0,0,0)^T.
\end{align*}
The calculation shows that $\frac{d}{d\lambda}\vert_{\lambda^j}\hat{\vel}^{j+1}$ is also a solution of a variational problem and independent from the Newton iteration $j$. In other words, $\Phi(\lambda)$ is linear and the Newton iteration should convergence in one step. For every timestep we therefore consider the following steps:  
\begin{enumerate}
    \item initial $\lambda^0=0$.   
    \item compute the derivative $\frac{d}{d\lambda}\Phi(\lambda)$ by finding a solution $(\hat{\delta\vel}^{n+1},\hat{\delta p}^{n+1},\hat{\delta\meanc}^{n+1},\hat{\delta\paramUpdate}^{n+1})$ of the variational problem 
        \begin{align*}
            L^n((\hat{\delta\vel}^{n+1},\hat{\delta p}^{n+1},\hat{\delta\meanc}^{n+1},\hat{\delta\paramUpdate}^{n+1}),(\velTest,q,h,\paramTest)) = (-\tau (\normal , \velTest),0,0,0)^T 
        \end{align*}    
        $\forall (\velTest,q,h,\paramUpdate)\in V_h(\surf_h^n)$.
    \item for every $j>0$ do the Newton iteration by
    \begin{enumerate}
        \item compute $\Phi(\lambda^j)=\int_{\surf_h^n} \hat{\vel}^{j+1}\cdot\normal$ by finding a solution   $(\hat{\vel}^{j+1},\hat{p}^{j+1},\hat{\meanc}^{j+1},\hat{\paramUpdate}^{j+1})$ of the variational problem 
        \begin{align*}
            L^n((\hat{\vel}^{j+1},\hat{p}^{j+1},\hat{\meanc}^{j+1},\hat{\paramUpdate}^{j+1}),(\velTest,q,h,\paramTest)) = r^n((\velTest,q,h,\paramTest)) + (-\tau \lambda^j(\normal , \phi),0,0,0)^T 
        \end{align*}
        $\forall (\velTest,q,h,\paramTest)\in V_h(\surf_h^n)$.
        \item compute the Newton update $\delta\lambda^{j+1} = -(\frac{d}{d\lambda}\Phi(\lambda^j))^{-1}\Phi(\lambda^j)$ and $\lambda^{j+1}=\lambda^j+\delta\lambda^{j+1}$.
        \item if $\int_\surf \hat{\vel}^{j+1}\cdot\normal \mathrm{d}\surf<\epsilon$ stop the Newton iteration, otherwise repeat with (a).    
    \end{enumerate}
    \item compute the new surface $\surf_h^{n+1}$ by updating its parametrization $\param^{n+1} = \param^n + \hat{\paramUpdate}^{j+1}$.  
    \item interpolate the velocity $\vel^{n+1}$ on $\surf^{n+1}$ by $\vel^{n+1}(x_k^{n+1})=\hat{\vel}^{j+1}(x_k^n)$ for $k=1,...,N_{Dof}$ where $N_{Dof}$ is the number of degrees of freedom of the triangulation. 
\end{enumerate}
This approach considers the volume constraint and the local inextensebility constraint only with respect to the old geometry $\surf_h^n$ and therefore does not guarantee $\DivP\vel^{n+1}=0$ and $\int_{\surf_h^{n+1}} \vel^{n+1}\cdot\normal \mathrm{d}\surf=0$. However, the convergence studies in the next section demonstrate the boundedness of these quantities and the applicability of the approach. 

\subsection{Remarks on the numerical approach}
As there are no theoretical results available for eqs. \eqref{eq:conservNavierStokes} - \eqref{eq:meshflow} or their discretization, the best to hope for are convergence properties known for sub-problems or to achieve better experimental convergence properties than previous numerical studies. Extending the results of \cite{hardering2021tangential} to the surface Stokes problem on stationary surfaces shows for a similar Taylor-Hood approach optimal order of convergence of degree $k+1$ for $\|\vel - \vel_h\|_{L^2}$ \cite{hardering2022Stokes}. This property and order $k$ convergence for $\|\DivP \vel_h\|_{L^2}$ have been shown numerically in \cite{Brandneretal_SIAMJSC_2022} for the surface Stokes problem on stationary surfaces. Results on evolving surfaces for these equations are not known. Shape relaxation by pure bending forces have been considered by different approaches. The closest to our scheme \cite{Dziuk_NM_2008} considers a piece-wise flat surface triangulation and leads to an experimental order of degree two for the considered Helfrich/Willmore energy without any constraints at the equilibrium state. A convergence proof for this result is still open. An extension of this approach to higher order (quadratic) surface approximations and volume and global area constraints in \cite{Bonitoetal_JCP_2010} does not provide any convergence results. Finally, the overall problem but without volume constraint has been considered in \cite{reuther2020numerical}. The discretization considers a piece-wise flat surface triangulation, a splitting in tangential and normal components, a Chorin-type projection for the tangential velocity with $k=1$ for velocity and pressure and artificial diffusion for the normal velocity. The considered examples show an experimental order of convergence of degree two in space and degree one in time for $\|\DivP \vel_h\|_{L^\infty(L^2)}$. In \cite{torres2019modelling} the full model is addressed. However, it is reformulated in a vorticity-stream function formulation, which allows only to deal with scalar-valued quantities but restricts the applicability to simply-connected surfaces \cite{Nitschkeetal_PAMM_2020}. Convergence studies are not provided in \cite{torres2019modelling} and can therefore not be considered for any comparison.

We will compare the properties of our discretization scheme with these known results. The described discretization scheme allows for straight forward parallelization, which is used in the following. The discretization is realized within the finite element toolbox AMDiS \cite{vey2007amdis,witkowski2015software} using the \textsc{Dune}-CurvedGrid library \cite{praetorius2020dunecurvedgrid}.  

\begin{figure}[h!]
    \centering
    \includegraphics[width=0.99\linewidth]{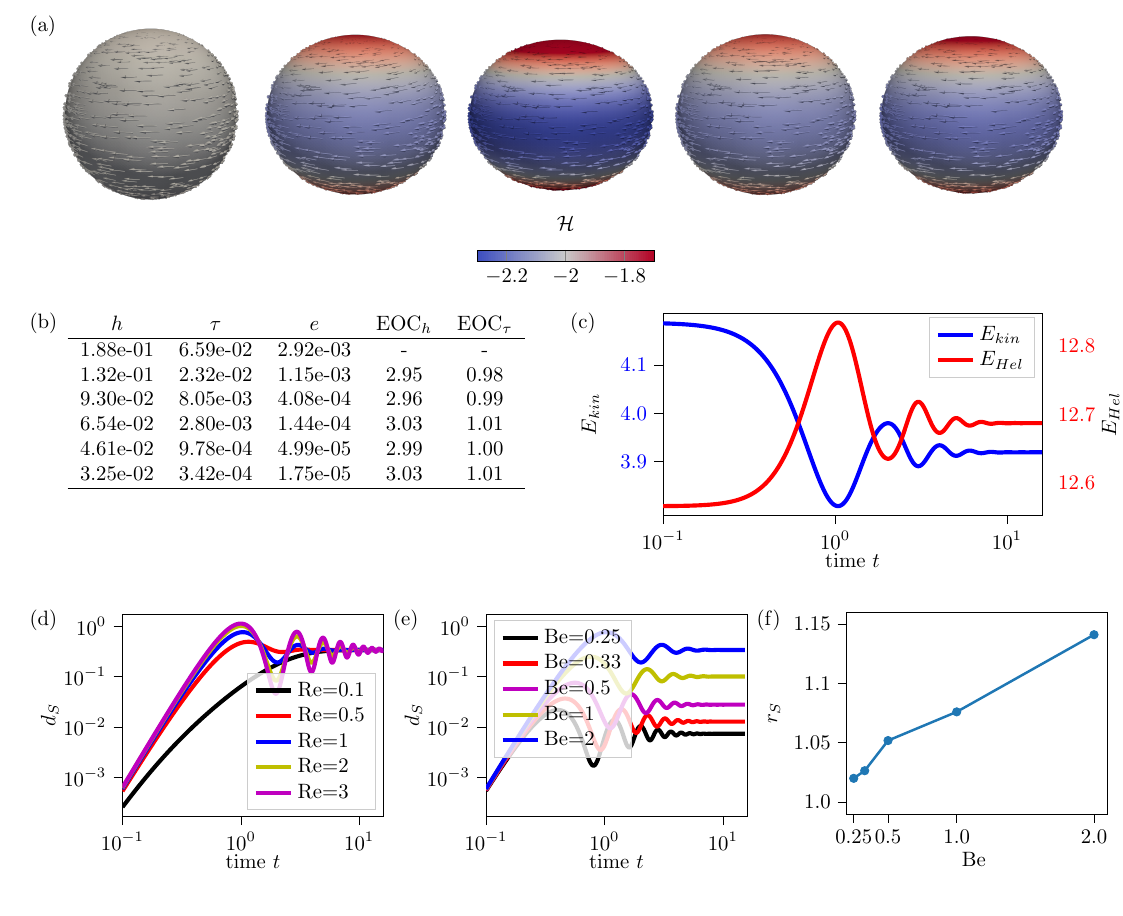}
    \caption{(a) Snapshots of the surface at $t=0,0.5,1,2,10$ for $\Re=1$ and $\Be=2$. Visualized by the glyphs is the tangential part of the surface velocity $\proj\vel$ where the length represents it's magnitude. By the color the mean curvature $\meanc$ is highlighted.  
    (b) Convergence study for the inextensibility error $e = \|\DivP \vel_h\|_{L^2}$ for the simulation shown in (a) with respect to $h$ and $\tau$. The results indicate $3$rd order in space and $1$st order in time. (c) Kinetic and Helfrich energy over time for the simulation  shown in (a). 
    (d,e) Derivation from the sphere $d_\surf = \int_\surf (\meanc-\meanc_{Ref})^2 \mathrm{d}\surf$ for different Reynolds numbers $\Re$ ($\Be=2$) and for difference bending capillary numbers $\Be$ ($\Re=1$), respectively. 
    (f) Relation between the ratio $r_\surf = \lim_{t\to \infty} (\max_{\mathbf{x}\in\surf} \Vert \mathbf{x} \Vert / \min_{\mathbf{x}\in\surf} \Vert \mathbf{x} \Vert)$ and $\Be$.}  
    \label{fig:killingField0}
\end{figure} 

\section{Results} 
\label{sec:results}

We consider three examples to demonstrate the applicability of the approach and the richness of the dynamics emerging from the solid-fluid duality: The first is a Killing vector field on a sphere where we drop the volume constraint. This example is taken from \cite{reuther2020numerical}, where it is used to demonstrate the interplay of tangential flow and surface deformation. The tangential velocity induces deformations towards ellipsoidal-like shapes. Due to the induced normal deformation energy dissipates which reduces the tangential velocity. With the numerical approach considered in \cite{reuther2020numerical}, which contains additional dissipation due to artificial diffusion, the shape relaxes back to a sphere with zero tangential velocity. 
However, it is already argued that a force balance between the bending forces of the Helfrich energy and the induced normal deformation due to tangential flow can be established. Our algorithm converges to these states which marks a qualitative difference to the previous results in \cite{reuther2020numerical}. The emerging shapes strongly depend on the bending capillary number $\Be$ and are validated by a semi-analytic approach assuming an axisymmetric solution. As in \cite{reuther2020numerical} we consider $\|\DivP \vel_h\|_{L^\infty(L^2)}$ as an appropriate error measure to consider convergence properties. The second example is also taken from \cite{reuther2020numerical} and considers a perturbed sphere with zero tangential velocity. The same setting is considered here, but with volume conservation. Detailed parameter studies in the Reynolds number $\Re$ and the bending capillary number $\Be$ show the enhanced convergence of the solution to its equilibrium for increasing $\Re$. The evolution is characterized by two significant shape changes. First, a biconcave-like shape with almost zero tangential velocity is approached. However, this is only a local minimum. Further evolution leads to a dumbbell shape and requires breaking symmetry. This second shape transition strongly depends on the flow properties and could not be achieved for the lowest considered $\Re$ within the used time frame. The biconcave and dumbbell shapes correspond to the axisymmetric equilibrium shapes emerging from oblate and prolate ellipsoids, respectively. These stationary shapes have been intensively discussed for the Helfrich energy neglecting any flow \cite{Seifert_AP_1997}. The parameters are chosen to highlight the influence of $\Re$ on the evolution. Again $\|\DivP \vel_h\|_{L^\infty(L^2)}$ is considered for establishing an experimental order of convergence. In addition we consider convergence properties in the Helfrich energy $E_H$ and in area and volume conservation. These convergence studies are done for $\Re = 1$ and $\Be=2$. We further demonstrate the influence of $\Re$ on mesh quality and area and volume conservation. The third example demonstrates the applicability of the approach to a more complex situation with extreme shape changes. We consider the Stanford bunny, which has been used for the incompressible surface Navier-Stoles model in \cite{Nitschkeetal_JFM_2012,Reutheretal_MMS_2015,Ledereretal_IJNME_2020} highlighting the connection between vortices and saddle points in the flow field and extrema in the Gaussian curvature of the surface. However, previous investigations are restricted to stationary surfaces. We here show the convergence to a biconcave shape and the induced tangential flow during the evolution. The complexity of the surface demonstrates the applicability of the considered mesh regularity approach but also points to its limits. For all example we use order $k=3$.

\subsection{Killing vector field}
\label{subsec:killing}

We consider the initial surface $\surf = \surf(0)$ to be a sphere with radius $R=1$, which corresponds to $|\Omega(0)| = 4.19$, and define a tangential rotating velocity field $\vel(0,\coord) = (\coordC_1,-\coordC_0,0)^T$ for $\coord=(\coordC_0,\coordC_1,\coordC_2)^T\in\surf$. With volume constraint this would be an equilibrium solution, a so-called Killing vector field, a velocity field which is not influenced by the rate of deformation tensor $\boldsymbol{\sigma}$. Such solutions have been numerically considered in \cite{Reutheretal_MMS_2015,Nitschkeetal_TFI_2017} and also recently analytically investigated \cite{Pruessetal_JEE_2021}. However, if the volume constraint is dropped shape deformations become possible and the tangential flow field induces a normal velocity which deforms the surface to an ellipsoid-like shape. This deformation leads to dissipation. We consider an initial mesh with $h=1.88e-1$. The simulation results for $\Re=1$ and $\Be=2$ are shown in Figure \ref{fig:killingField0}(a). The velocity $\vel$ is visualized by arrows, showing the direction and magnitude and the shape is highlighted by plotting the mean curvature $\meanc$. The surface oscillates around the stable ellipsoidal-like shape until a balance between the acting forces is reached. We additionally plot the kinetic energy $E_{kin}=\frac{1}{2\Re}\int_\surf \Vert \vel \Vert^2 \mathrm{d}\surf$ and the Helfrich energy $E_H$ over time. As long as the surface is deforming the rate of deformation tensor acts on the normal velocity and the kinetic energy dissipates.

For the parameters $\Re=1$ and $\Be=2$ a convergence study with respect to the grid size $h$ and the timesteps width $\tau$ is shown in Figure \ref{fig:killingField0}(b), where $\tau\sim h^3$. We consider the inextensibility error $e = \|\DivP \vel_h\|_{L^\infty(L^2)}$. We obtain third order convergence with respect to $h$ and first order convergence with respect to $\tau$. For the spatial error this corresponds to the experimental order of convergence obtained in \cite{Brandneretal_SIAMJSC_2022} for the surface Stokes model on stationary surfaces, which is expected to be optimal for this problem.

In Figure \ref{fig:killingField0}(d) and (e) the behavior for different $\Re$ and $\Be$ is shown. We measure $d_{\surf} = \int_\surf (\meanc-\meanc_{Ref})^2 \mathrm{d}\surf$ the derivation from the sphere where $\meanc_{Ref}$ is the mean curvature of a sphere with the same area as the surface $\surf$. Similar to \cite{reuther2020numerical} we observe larger oscillations for larger Reynolds numbers $\Re$ and smaller bending capillary numbers $\Be$. Decreasing $\Re$ dampens the oscillations amplitude but has little influence on the frequency, except for $\Re=0.1$, for which the damping is so high that no oscillations occur. $\Re$ has no influence on the equilibrium state. It is solely determined by $\Be$ and the initial conditions. As smaller $\Be$ as more deformed the equilibrium shape gets. Also the frequency of the oscillations increases with decreasing $\Be$. We characterize the equilibrium shape by the aspect ratio $r_\surf = \lim_{t\to \infty} (\max_{\mathbf{x}\in\surf} \Vert \mathbf{x} \Vert / \min_{\mathbf{x}\in\surf} \Vert \mathbf{x} \Vert)$. Its dependency on $\Be$ is plotted in Figure \ref{fig:killingField0}(f) and quantifies the divergence from a sphere ($r_\surf = 1.0$) for increasing $\Be$. 

Assuming axisymmetric solutions, characterizing the tangential flow field as a Killing vector field and assuming a stationary ellipsoidal geometry $\surf_{eq}$, which is determined by the force balance between normal forces induced by the tangential flow and bending forces allows to simplify the problem setting. Without loss of generality we consider the $z$-axis as the symmetry axis. The equilibrium velocity field is define by a rigid body rotation
\begin{align}
    \label{eq:errorKillingU}
    \vel &= \omega \boldsymbol{e}_z\times\coord = \omega (\coordC_1,-\coordC_0,0)^T
\end{align}
where $\omega\in\R$ is the angular velocity and $\boldsymbol{e}_z$ the unit vector of the $z$-axis. The velocity $\vel$ fulfills $\DivP \vel = 0$ and has no influence on the rate-of-deformation tensor $\stress$. The induced centrifugal forces $\boldsymbol{F}_z\coloneqq\nabla_{\vel}\vel=-\omega^2\boldsymbol{e}_z \times (\boldsymbol{e}_z \times \coord)$ are balanced by surface pressure and bending forces. This simplifies eqs. \eqref{eq:navierstokes} which can be separated in a tangential and a normal part by multiplying with $\proj$ and $\normal$, respectively. We obtain
\begin{align}
    \omega^2\proj \boldsymbol{F}_z &= -\GradSurf p  & \text{ on } \surf_{eq}     \label{eq:errorKillingP} \\
    \omega^2\boldsymbol{F}_z\cdot\normal &= - p\meanc + \frac{1}{\Be}(-\Delta_\surf\meanc -\meanc(\Vert\shop\Vert^2 - \frac{1}{2}\operatorname{tr}(\shop)^2) & \text{ on }\surf_{eq}. \label{eq:errorKillingH}
\end{align}

For a given angular velocity $\omega$ we consider $\surf_{eq}$ and $p$ as unknowns. Transforming eqs. \eqref{eq:errorKillingP} and \eqref{eq:errorKillingH} in spherical coordinates leads to a system of ordinary differential equations, where the radius $R$ and the pressure $p$ remain as unknowns. However, the equations are highly nonlinear, contain derivatives of the radius up to fourth order and singularities at the poles. We were unable to solve them, either analytically nor numerically. However, we can still consider the equations above to verify the simulation results shown in \ref{subsec:killing}. We consider the discrete solution $(\vel_h,p_h,\param_h,\meanc_h)$ at the final time step and define the following residuals for a given $\omega$ 
\begin{align*}
    r_0 &= \Vert \vel_h\cdot \normal_h \Vert_{L^2(\surf_h)} \\
    r_1 &= \Vert \vel_h-\omega \boldsymbol{e}_z\times\coord \Vert_{L^2(\surf_h)} \\
    r_2 &= \Vert \omega^2\proj \boldsymbol{F}_z + \GradSurf p_h \Vert_{L^2(\surf_h)} \\
    r_3 &= \Vert \omega^2\boldsymbol{F}_z\cdot\normal + p_h\meanc_h - \frac{1}{\Be}(-\Delta_\surf\meanc_h -\meanc_h(\Vert\shop_h\Vert^2 - \frac{1}{2}\operatorname{tr}(\shop_h)^2) \Vert_{L^2(\surf_h)}. 
\end{align*}
They correspond to the normal component of the velocity, the tangential component of the velocity, eq. \eqref{eq:errorKillingP} and eq. \eqref{eq:errorKillingH}, respectively. If these quantities are small our numerical solution approximates the semi-analytic axisymmetric solution on the stationary ellipsoidel geometry $\surf_{eq}$. The results are shown in Table \ref{tab:residuals} for different $\Re$ and $\Be$. The residuals are sufficiently small to claim approximation of the semi-analytic solution which confirms the validity of the numerical scheme. Due to the different values obtained for $\omega$ the results also confirm our finding on the dependency of the equilibrium shape on $\Be$. The results also quantitatively show that $\omega$ and thus also the equilibrium shape, do not depend on $\Re$. The values remain constant for all considered $\Re$.
\begin{table}[h!]
    \centering
    \begin{tabular}{ccccccc}
     $\Re$&$\Be$ &$\omega$  & $r_0$    & $r_1$    & $r_2$    & $r_3$ \\ \hline
     $1.0$&$0.25$& $0.9894$ & $3.4\float{5}$ & $1.9\float{4}$ & $1.9\float{3}$ & $1.6\float{2}$ \\
     $1.0$&$0.33$& $0.9862$ & $3.4\float{5}$ & $6.7\float{4}$ & $1.9\float{3}$ & $1.2\float{2}$ \\
     $1.0$&$0.5$ & $0.9801$ & $3.3\float{5}$ & $1.1\float{3}$ & $2.0\float{3}$ & $8.3\float{3}$ \\
     $1.0$&$1.0$ & $0.9629$ & $3.2\float{5}$ & $6.5\float{4}$ & $1.9\float{3}$ & $4.6\float{3}$ \\
     $1.0$&$2.0$ & $0.9356$ & $3.3\float{5}$ & $6.4\float{4}$ & $1.9\float{3}$ & $2.8\float{3}$ \\ 
     $0.1$&$2.0$ & $0.9356$ & $5.1\float{5}$ & $1.1\float{4}$ & $2.0\float{3}$ & $2.8\float{3}$ \\
     $0.5$&$2.0$ & $0.9356$ & $3.5\float{5}$ & $7.5\float{4}$ & $1.9\float{3}$ & $2.7\float{3}$ \\     
     $2.0$&$2.0$ & $0.9356$ & $2.7\float{3}$ & $3.5\float{3}$ & $3.4\float{3}$ & $3.3\float{3}$ \\
     $3.0$&$2.0$ & $0.9356$ & $1.4\float{2}$ & $1.8\float{2}$ & $1.1\float{2}$ & $5.7\float{3}$ \\     
    \end{tabular}
    \caption{The angular velocity $\omega$ and residual errors for a variation of $\Re$ and $\Be$.}
    \label{tab:residuals}
\end{table}

\begin{figure}[ht!]
    \centering
    \includegraphics[width=0.99\linewidth]{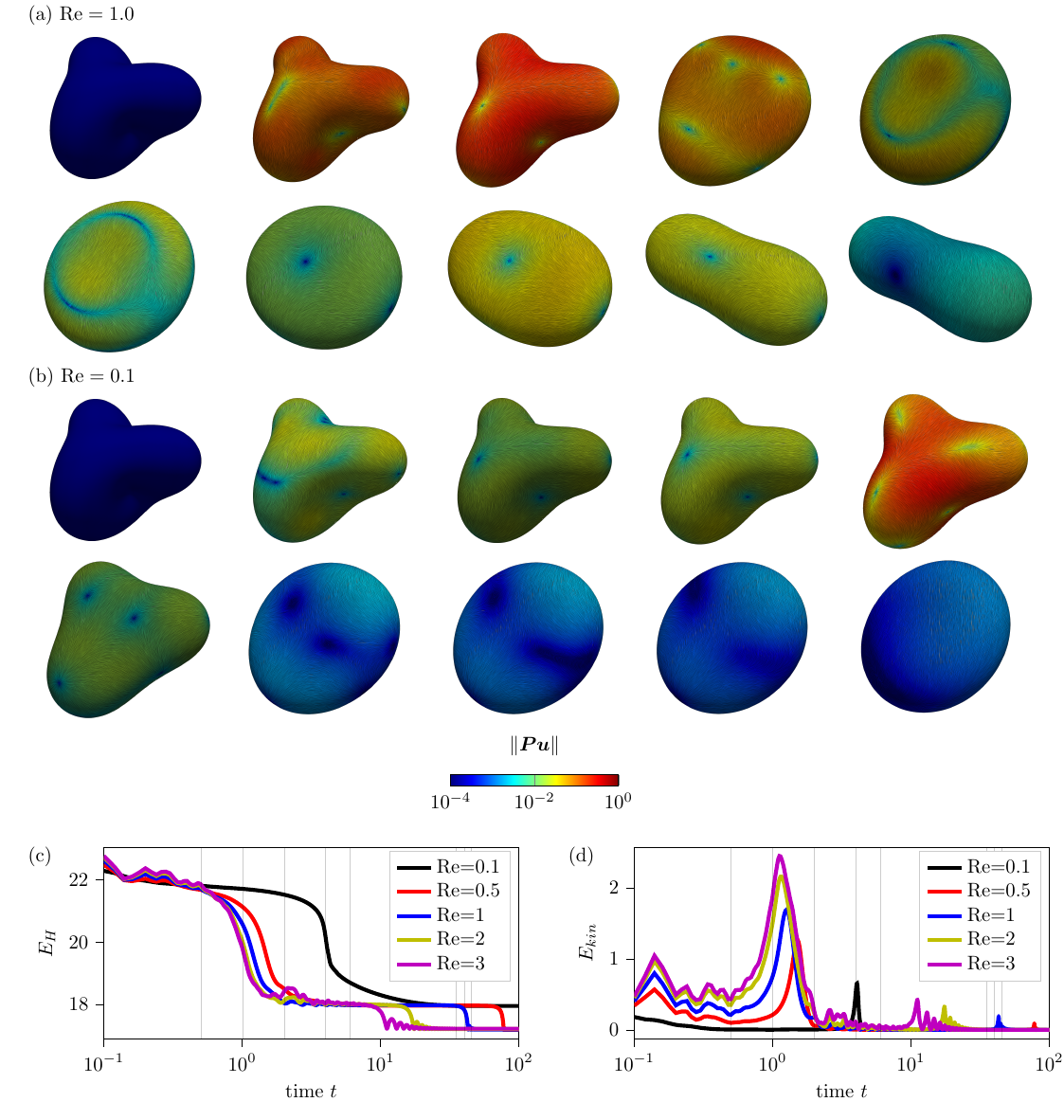} 
    \caption{Relaxation of a perturbed sphere for Reynolds numbers $\Re=1$ (a) and $\Re=0.1$ (b) at snapshots at the time $t=0.0,0.5,1,2,4,6,35,40,45,100$. Visualized is the tangential flow field using LIC with the magnitude color coded with logarithmic scale. (c) shows the Helfrich energy $E_H$ over time for different $\Re$. (d) shows the kinetic energy $E_{kin}$ over time for different $\Re$. The time instances considered in (a) and (b) are marked in the plots.} 
    \label{fig:conservative1}
\end{figure}

\subsection{Relaxation of a perturbed sphere}

We consider a perturbed sphere $\surf = \surf(0)$ with $\vel =0$ as an initial condition. The perturbation is given by the parametrization $\param(\theta,\psi) = r(\theta,\psi)\param_S(\theta,\psi)$ where $\param_S(\theta,\psi)$ is the parametrization of the unit sphere and $r(\theta,\psi) = 1 + r_0\cos{\theta}\sin{3\psi}$ a space dependent radius with $r_0=0.4$. The considered perturbed sphere has volume $|\Omega| = 4.70$ and area $|\surf| = 15.64$. As a result of the bending forces, the shape relaxes and thereby induces flow. We consider again an initial mesh with $h=1.88e-1$. Figure \ref{fig:conservative1}(a) and (b) show the simulation results for $\Re=1$ and $\Re=0.1$, respectively. The bending capillary number is $\Be=2$. 
We consider a LIC visualization with color coding showing the magnitude of the tangential part of the velocity $\|\proj \vel\|$. This highlights the correlation between shape deformation and tangential flow. The shown time instances are equal for both simulations and chosen to highlight strong shape deformations during the evolution. The shape first relaxes towards a biconcave shape with almost zero tangential velocity. Further relaxation leads to a change in geometry and convergence towards a dumbbell shape. The tangential velocity in the final state is zero. This equilibrium shape is in agreement with the phase diagram discussed for the Helfrich energy neglecting any flow in \cite{Seifert_AP_1997}. The surface area and enclosed volume in the considered example lead to a reduced volume of $0.82$, see \cite{Seifert_AP_1997}, and thus correspond to an equilibrium dumbbell shape in the phase diagram. 

However, the transition to the dumbbell shape is only achieved for $\Re = 1$ and not observed within the considered time frame for $\Re = 0.1$. In order to explain this behaviour we plot the Helfrich energy $E_H$ and the kinetic energy $E_{kin}$ over time for different $\Re$ in Figure \ref{fig:conservative1}(c) and (d), respectively. The different time instances are marked in the plots. The transition from a biconcave to a dumbbell shape is associated with a steep decrease in $E_H$. With lower $\Re$, this transition occurs later in time. For larger $\Re$ we also observe shape oscillations after such drastic shape changes. However, any change in $E_H$ is also associated with an increase in $E_{kin}$. We conclude that surface flow enhances shape relaxation. In the considered example it even helps to escape the local minimum of a biconcave shape.  

\begin{figure}[h!]
    \centering
    \includegraphics[width=0.99\linewidth]{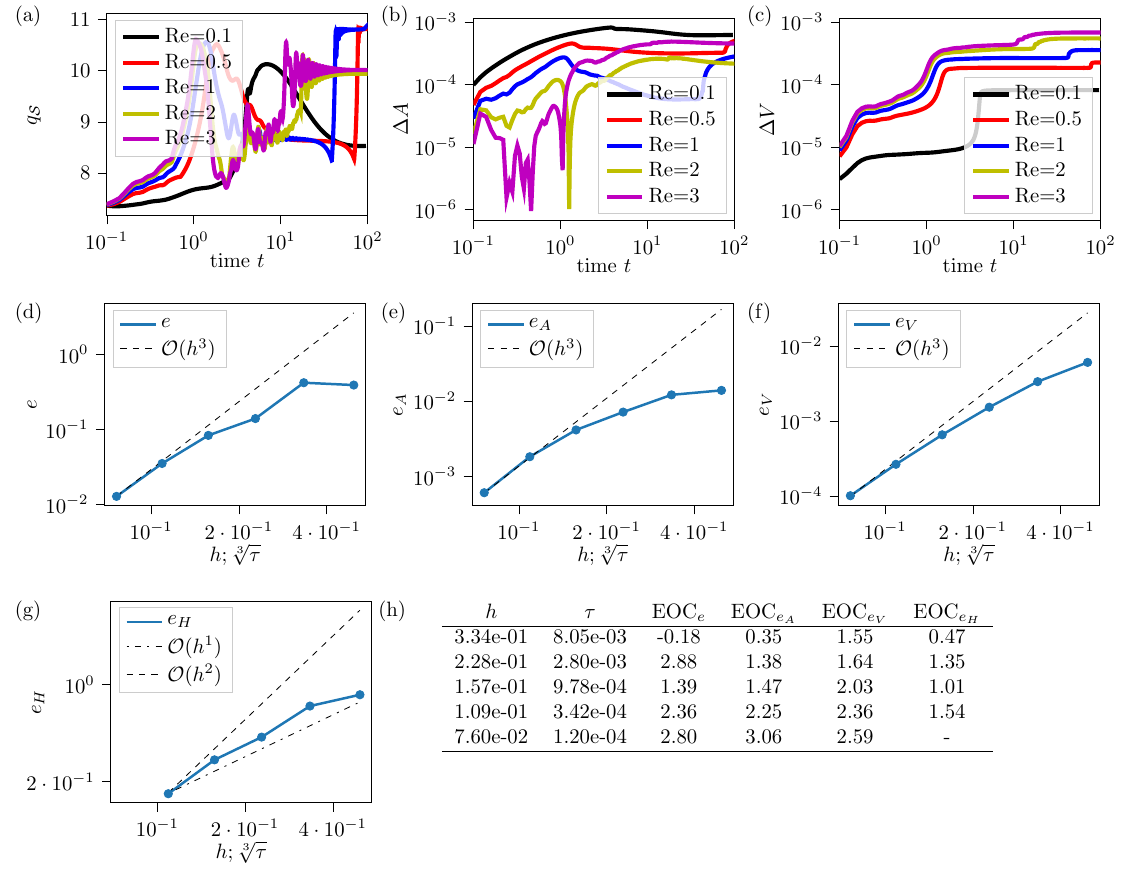}
    \caption{(a), (b) and (c) show the quasi uniformity of the mesh, the area conservation error and the volume conservation error over time for different $\Re$, respectively. For $\Re = 1$ in (d) the inextensibility error $e = \|\DivP \vel_h\|_{L^{\infty}(L^2)}$, in (e) the error in area conservation $e_A = \| \Delta A \|_{L^{\infty}}$, in (f) the error in volume conservation $e_V = \| \Delta V \|_{L^\infty}$, and in (g) the error of the Helfrich energy $e_H = \Vert E_H \Vert_{L^\infty}$ with respect to $h$ and $\tau$ is shown. (h) shows the experimental convergence orders of the errors of (d)-(g).}  
    \label{fig:conservative2}
\end{figure}

We next use these simulations to demonstrate the numerical properties of the algorithm. This includes the quality of the mesh, as well as the global constraints on area and volume. Figure \ref{fig:conservative2}(a), (b) and (c) show the quasi uniformity of the mesh, which is measured by $q_{\surf} = (\max_{T\in\mathcal{T}} h_T) / (\min_{T\in\mathcal{T}} h_T)$, the error in area conservation $\Delta A(t)= \vert \vert \surf(t) \vert - \vert \surf(0)\vert \vert / \vert \surf(0) \vert$ and the error in volume conservation $\Delta V(t)= \vert \vert \Omega(t) \vert - \vert \Omega(0)\vert \vert / \vert \Omega(0) \vert$. All these quantities stay bounded within the considered time frame. 
The volume is thereby computed as $\vert\Omega\vert = \frac{1}{3} \int_\Omega \Div(\coord) \mathrm{d}\coord = \frac{1}{3}\int_\surf \coord\cdot \normal \mathrm{d}\surf$, see \cite{hughes1996application,bao2021structure}. For $\Re=1$ and $\Be=2$ a convergence study is shown in Figure \ref{fig:conservative2}(d), (e) and (f). As before we reduced the time step width with the grid size by $\tau\sim h^3$. The inextensebility error $e$ and the $L^\infty$-norms of the area and volume error $e_A$ and $e_V$ are measured. As before we observe third order convergence for $e$ with respect to $h$ and first order with respect to $\tau$. This result suggest that the volume constraint does not affect the convergence properties.

For $e_A$ and $e_V$ the results also indicate an approaching third order convergence with respect to $h$. However, this is only achieved as long as the error bound of the Newton method is smaller than the interpolation error. In the current situation the parameter $\epsilon$ in the Newton method is chosen as $\epsilon=1\float{6}$. Figure \ref{fig:conservative2}(g) shows in addition convergence properties for $E_H$, again considered in the $L^\infty$-norm. These results are obtained with respect to the numerical solution of the most refined mesh. The corresponding values and the computed experimental orders of convergence (EOC) are shown in Figure \ref{fig:conservative2}(h). The obtained order of convergences for $E_H$ is between one and two. For pure Willmore flow without any constraints and only considering the equilibrium state \cite{Dziuk_NM_2008} showed second order convergence. Compared with this experimental result for a (still complex but) strongly simplified sub-problem these convergence results also confirm the validity of the proposed approach.

\begin{figure}[h!]
    \centering
    \includegraphics[width=\linewidth]{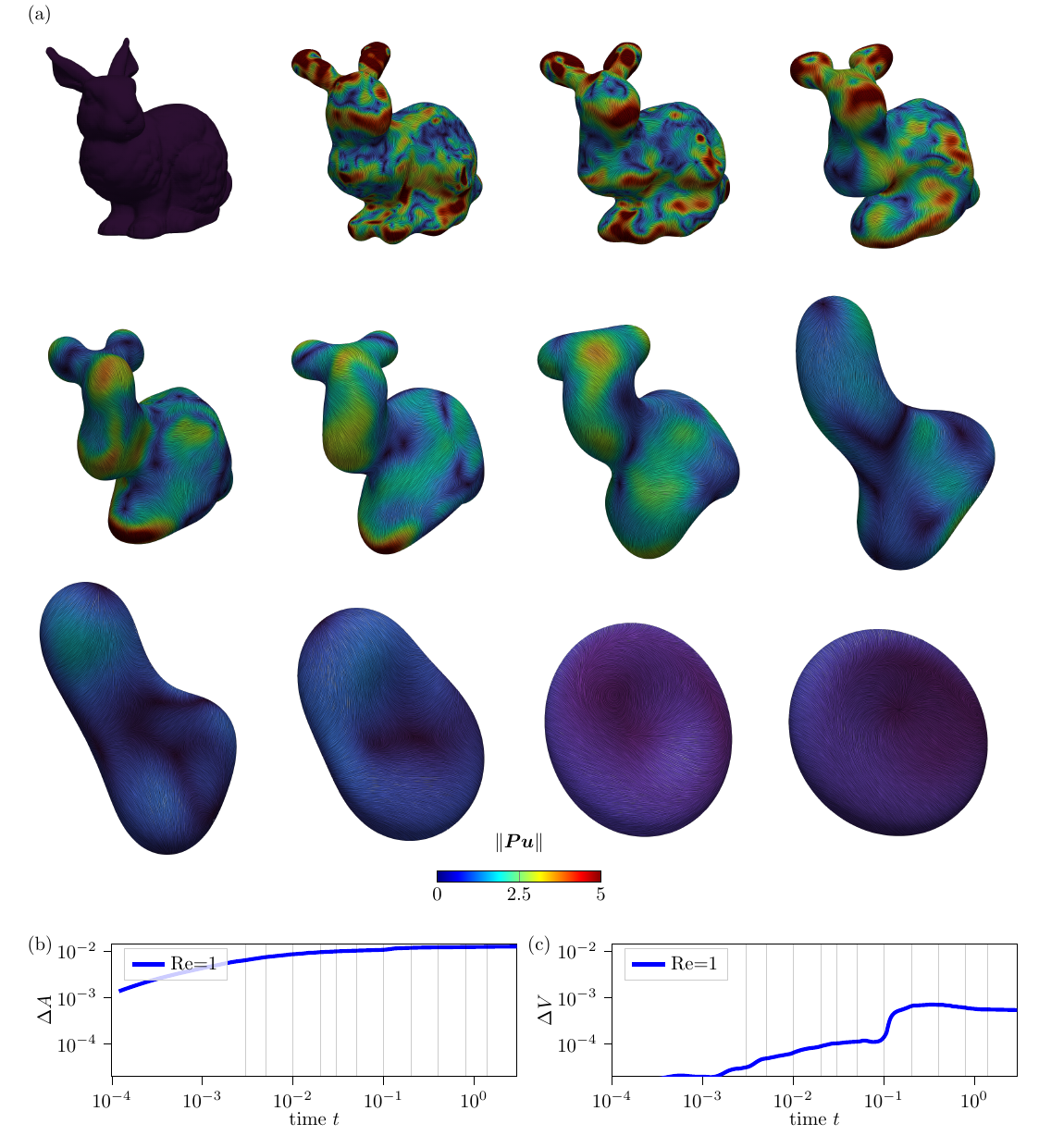}
    \caption{(a) Relaxation of the Stanford bunny for Reynolds numbers $\Re=1$ and $\Be=0.5$ at snapshots at the time $t=0.0,  0.003, 0.005, 0.01, 
         0.02 ,0.03 , 0.05 , 0.1, 
         0.2 , 0.4 ,  0.8 , 1.4$. Visualized is the flow generated by the tangential part of the surface $\proj\vel$ by a LIC filter. A corresponding movie is provided in the Supplement.
    (b) and (c) show the area conservation error and the volume conservation error over time. The time instances considered in (a) are marked in the plots.} 
    \label{fig:bunny1}
\end{figure}

\subsection{Relaxation of the Stanford bunny}

Let $\surf = \surf(0)$ be the Stanford bunny with initial velocity $\vel=0$. The Stanford bunny has volume $|\Omega| = 0.146$ and area $|\surf| = 1.95$ and is triangulated by $70.118$ elements which leads to around $2.5$ million degrees of freedom in the discrete system. The initial mesh size is $h=2\float{2}$ and the corresponding time step is chosen as $\tau=1\float{5}$. As for the perturbed sphere we consider $\Re=1$ and $\Be=2$. 

The Stanford bunny has previously been used to demonstrate the validity of numerical algorithms to solve the incompressible surface Navier-Stoles model, see \cite{Nitschkeetal_JFM_2012,Reutheretal_MMS_2015,Ledereretal_IJNME_2020}. However, in these simulations the surfaces were stationary and only a tangential flow was considered. The results for the full model are shown in Figure \ref{fig:bunny1}(a).

Driven by the Helfrich energy the shape relaxes and a tangential flow is induced. This complex solid-fluid duality together with area and volume conservative leads to a complex evolution of the Stanford bunny which finally converges to a biconcave shape with zero tangential velocity. Comparing the equilibrium shape with the phase diagram in \cite{Seifert_AP_1997}, shows for the considered volume and area a reduced volume of $0.57$ and is therefor below the transition to the dumbbell shape and within the transition zone between a biconcave and a stromatocyte shape. Figure \ref{fig:bunny1}(b) and (c) show the error in area and volume conservation, respectively. Again, the shown time instances are marked in these plots. While larger in magnitude than in the previous examples, the error bounds are still acceptable for such complex shape evolution with extreme curvatures and large velocity differences. 

\begin{figure}[h!]
    \centering
    \includegraphics[width=\linewidth]{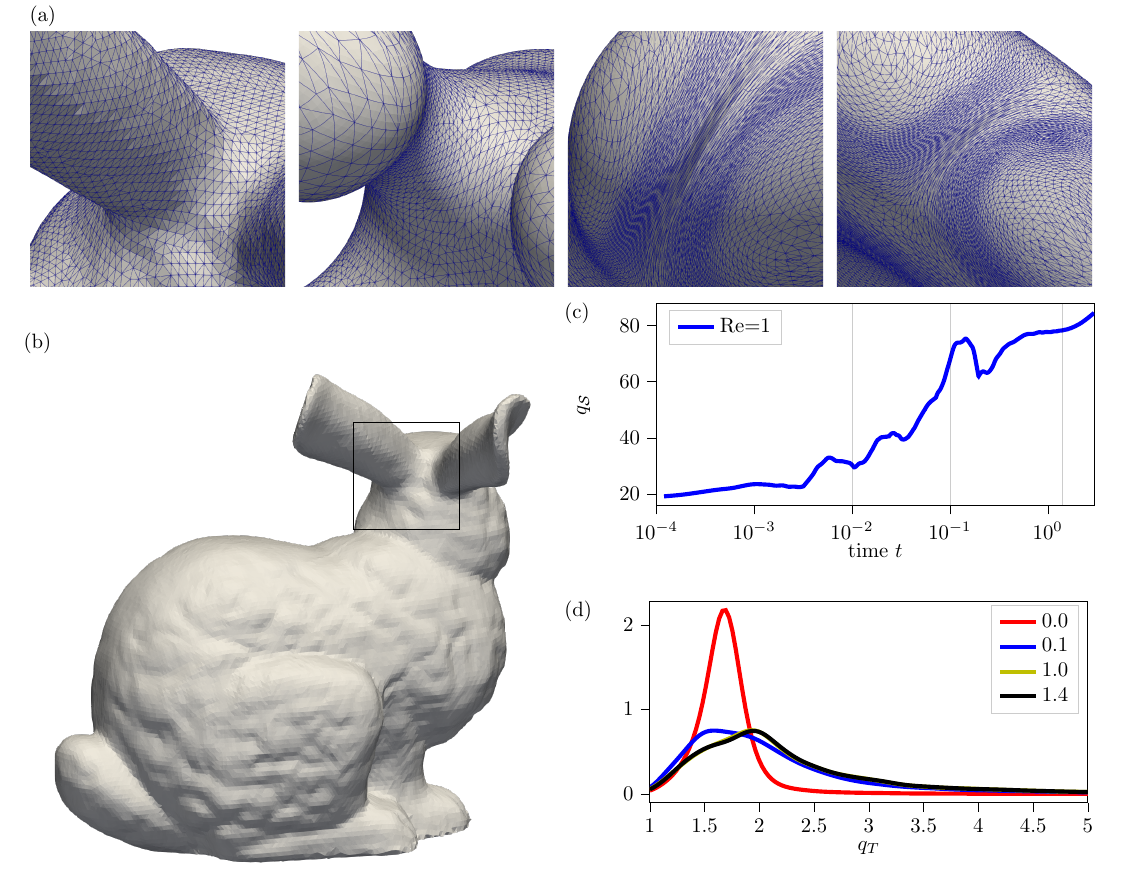}
    \caption{(a) The mesh of the simulation shown in figure \ref{fig:bunny1} (a) at the time $t=0.0, 0.01, 0.1, 1.4$. The Snapshots are follows on the region between the ears of the bunny highlighted in (b). (c) shows the quasi uniformity of the mesh. The time instances considered in (a) are marked in the plot. (d) discrete probability density of the triangle quality $q_T$ for all triangles of the mesh and for the time instances shown in (a).} 
    \label{fig:bunny2}
\end{figure}

This complex evolution also provides an ideal test case for the considered local mesh redistribution algorithm based on \cite{barrett2008parametric}. We therefore plot a region of the bunny and the associated triangulation, which is strongly deformed during evolution as shown in Figure \ref{fig:bunny2}(a) and (b). The initial mesh size of $h=2\float{2}$ increases in this region to a maximum of $1\float{1}$. The differences in mesh size are clearly visible in Figure \ref{fig:bunny2}(a). This is quantified in Figure \ref{fig:bunny2}(c) showing the quasi uniformity of the mesh, which increases over time and reaches much larger values than in previous examples. While not optimal with respect to efficiency, these differences in size do not influence the performance of the linear solver. A second criteria is the regularity of the triangles. 
We therefore consider the edge length of each surface triangle and compute $q_T = (\max_{F,F\cap T \neq \emptyset} \vert F \vert)/(\min_{F,F\cap T \neq \emptyset} \vert F \vert) $ where $F$ denotes the edges of the triangulation $\mathcal{T}$. An equilateral triangle provides a lower bound with $q_T=1$. Figure \ref{fig:bunny2}(d) shows a histogram for $q_T$ for the considered time instances shown in 
Figure \ref{fig:bunny2}(a). While a clear shift towards more distorted triangles can be seen, there are no extreme deformations which again validates the proposed algorithm.

\section{Discussion}
\label{sec:Discussion}

We have considered a model for fluid deformable surfaces. The model shows a solid-fluid duality and combines an incompressible surface Navier-Stokes equation to model the fluid properties in tangential and normal direction with a forcing term resulting from a Helfrich energy to consider bending forces in normal direction and restricts the evolution by enforcing a constant enclosed volume. The system of highly nonlinear partial differential equations is formulated in an Eulerian-Lagrangian approach and is numerically solved by surface finite elements (SFEM) which extends previous approaches for surface Stokes equations \cite{Brandneretal_SIAMJSC_2022} with algorithms for Willmore flow \cite{Dziuk_NM_2008} and a Lagrange multiplier approach to enforce the volume constraint \cite{torres2019modelling}. In order to deal with strong shape deformations a mesh redistribution in the tangential direction is enforced using ideas of \cite{barrett2008parametric}. A key to the success of the approach is the use of higher order surface approximations within the \textsc{Dune}-CurvedGrid library \cite{praetorius2020dunecurvedgrid}, which allows to optain higher order approximations of geometric quantities and to use Taylor-Hood elements for the distretization of the surface Navier-Stokes model. Having this flexibility allows us to account for the first numerical analysis results for surface vector-Laplacians on stationary surfaces \cite{Olshanskiietal_SIAMJSC_2018,Hansboetal_IMAJNA_2020,hardering2021tangential} which suggest the necessity of higher order normal approximations to ensure optimal convergence. 
We numerically showed various solution properties and demonstrated the applicability of the numerical scheme by showing experimental orders of convergence which have been analytically shown or numerically obtained for various sub-problems and can be assumed to be optimal. These results indicate that optimal orders for stationary surfaces can also be realized on evolving surfaces and ask for analytical justification. Some results in this direction, which extend the results of \cite{demlow2006localized} for scalar valued surface partial differential equations to evolving surfaces with prescribed evolution or simple evolution laws, exist, see \cite{Kovacs_IAMJNA_2018} and \cite{Elliottetal_NM_2022}, respectively. 
However, numerical analysis results for the vector valued surface partial differential equations are open. One drawback of the proposed algorithm is the first order convergence in time which leads to severe time step restrictions. The challenges to overcome this are highly complex. Some approaches exist for scalar-valued surface partial differential equations on evolving surfaces, see \cite{DLM_IMAJNA_2012,LMV_IMAJNA_2013,Kovacs_IAMJNA_2018} and also for geometric evolution laws higher order time discretizations have been proposed, see \cite{BR_IFB_2012,Kovacsetal_NM_2019}. However, extending these results to the considered problem of fluid deformable surfaces is a research topic on its own, which has not been addressed.

The proposed numerical approach provides the basis for reliable numerical investigations of the rich applications for fluid deformable surfaces in biology. This requires extension of the model to include active components and liquid crystal properties, see \cite{Boisetal_PRL_2011,Salbreuxetal_PRE_2017,Mickelinetal_PRL_2019,Apaza_SM_2018,Nitschkeetal_PRF_2019,Pearceetal_PRL_2019,Mietkeetal_PNAS_2019,Nitschkeetal_PRSA_2020,Ranketal_PF_2021,Al-Izzietal_SCDB_2021,Hoffmannetal_SA_2022} for attempts in this direction. However, all of them only consider special situations: either stationary shapes, no hydrodynamics, axisymmetric settings, questionable numerical approaches or no simulation results at all.  \\

\textbf{Acknowledgements.} This work was supported by the German Research Foundation (DFG) through grant VO 899/22 within the Research Unit ”Vector- and Tensor-Valued Surface PDEs” (FOR 3013). We further acknowledge computing resources provided by ZIH at TU Dresden and by JSC at FZ J\"ulich, within projects WIR and PFAMDIS, respectively.

%\appendix 
%\section{Equation for the equilibrium state}
%\label{sec:equilibrium}

\bibliographystyle{elsarticle-num}
\bibliography{lib} 

\begin{thebibliography}{10}
\expandafter\ifx\csname url\endcsname\relax
  \def\url#1{\texttt{#1}}\fi
\expandafter\ifx\csname urlprefix\endcsname\relax\def\urlprefix{URL }\fi
\expandafter\ifx\csname href\endcsname\relax
  \def\href#1#2{#2} \def\path#1{#1}\fi

\bibitem{torres2019modelling}
A.~Torres-S{\'a}nchez, D.~Mill{\'a}n, M.~Arroyo, Modelling fluid deformable
  surfaces with an emphasis on biological interfaces, J. Fluid Mech. 872 (2019)
  218--271.

\bibitem{Voigt_JFM_2019}
A.~Voigt, Fluid deformable surfaces, J. Fluid Mech. 878 (2019) 1--4.

\bibitem{reuther2020numerical}
S.~Reuther, I.~Nitschke, A.~Voigt, A numerical approach for fluid deformable
  surfaces, J. Fluid Mech. 900 (2020) R8.

\bibitem{Mayeretal_Nature_2010}
M.~Mayer, M.~Depken, J.~Bois, F.~J\"ulicher, S.~Grill, {Anisotropies in
  cortical tension reveal the physical basis of polarizing cortical flows},
  {Nature} {467} ({2010}) {617--621}.

\bibitem{Heisenbergetal_Cell_2013}
C.-P. Heisenberg, Y.~Bellaiche, {Forces in tissue morphogenesis and
  patterning}, {Cell} {153} ({2013}) {948--962}.

\bibitem{Arroyoetal_PRE_2009}
M.~Arroyo, A.~DeSimone, Relaxation dynamics of fluid membranes, Phys. Rev. E 79
  (2009) 031915.

\bibitem{Salbreuxetal_PRE_2017}
G.~Salbreux, F.~J\"ulicher, Mechanics of active surfaces, Phys. Rev. E 96
  (2017) 032404.

\bibitem{Nitschkeetal_JFM_2012}
I.~Nitschke, A.~Voigt, J.~Wensch, A finite element approach to incompressible
  two-phase flow on manifolds, J. Fluid Mech. 708 (2012) 418--438.

\bibitem{Barrettetal_PRE_2015}
J.~W. Barrett, H.~Garcke, R.~N\"urnberg, Numerical computations of the dynamics
  of fluidic membranes and vesicles, Phys. Rev. E 92 (2015) 052704.

\bibitem{Saffmanetal_PNAS_1975}
P.~G. Saffman, M.~Delbr\"uck, Brownian motion in biological membranes, Proc.
  Nat. Acad. Sci. (USA) 72 (1975) 3111--3113.

\bibitem{Mietkeetal_PNAS_2019}
A.~Mietke, F.~J\"ulicher, I.~F. Sbalzarini, Self-organized shape dynamics of
  active surfaces, Proc. Nat. Acad. Sci. (USA) 116 (2019) 29--34.

\bibitem{Sahuetal_JCP_2020}
A.~Sahu, Y.~Omar, R.~Sauer, K.~Mandadapu, {Arbitrary Lagrangian-Eulerian finite
  element method for curved and deforming surfaces. I. General theory and
  application to fluid interfaces}, J. Comput. Phys. 407 (2020) 109253.

\bibitem{Grossetal_SIAMJNA_2018}
S.~Gross, T.~Jankuhn, A.~A. Olshanskii, A.~Reusken, A trace finite element
  method for vector-{L}aplacians on surfaces, SIAM J. Numer. Anal. 56 (2018)
  2406--2429.

\bibitem{Hansboetal_IMAJNA_2020}
P.~Hansbo, M.~G. Larson, K.~Larsson, Analysis of finite element methods for
  vector {L}aplacians on surfaces, IAM J. Numer. Anal. 40 (2020) 1652--1701.

\bibitem{hardering2021tangential}
H.~Hardering, S.~Praetorius, Tangential errors of tensor surface finite
  elements, IMA J. Numer. Anal. (2021) drag015.

\bibitem{Olshanskiietal_SIAMJSC_2018}
A.~A. Olshanskii, A.~Quaini, A.~Reusken, V.~Yushutin, A finite element method
  for the surface {Stokes} problem, SIAM J. Sci. Comput. 40 (2018)
  A2492--A2518.

\bibitem{hardering2022Stokes}
S.~Praetorius, personal communication.

\bibitem{Fries_IJNMF_2018}
T.-P. Fries, Higher-order surface {FEM} for incompressible {Navier–Stokes}
  flows on manifolds, Int. J. Numer. Meth. Fluids 88 (2018) 55--78.

\bibitem{Reutheretal_PF_2018}
S.~Reuther, A.~Voigt, Solving the incompressible surface {Navier–Stokes}
  equation by surface finite elements, Phys. Fluids 30 (2018) 012107.

\bibitem{Grossetal_JCP_2018}
B.~Gross, P.~J. Atzberger, Hydrodynamic flows on curved surfaces: {S}pectral
  numerical methods for radial manifold shapes, J. Comput. Phys. 371 (2018)
  663--689.

\bibitem{nestler2019finite}
M.~Nestler, I.~Nitschke, A.~Voigt, A finite element approach for vector- and
  tensor-valued surface {PDEs}, J. Comput. Phys. 389 (2019) 48--61.

\bibitem{Ledereretal_IJNME_2020}
P.~L. Lederer, C.~Lehrenfeld, J.~Sch\"{o}berl, Divergence-free tangential
  finite element methods for incompressible flows on surfaces, Int. J. Numer.
  Meth. Eng. 121 (2020) 2503--2533.

\bibitem{Ranketal_PF_2021}
M.~Rank, A.~Voigt, Active flows on curved surfaces, Phys. Fluids 33 (2021)
  072110.

\bibitem{Brandneretal_SIAMJSC_2022}
P.~Brandner, T.~Jankuhn, S.~Praetorius, A.~Reusken, A.~Voigt, Finite element
  discretization methods for velocity-pressure and stream function formulations
  of surface {S}tokes equations, SIAM J. Sci. Comput. 44 (2022) A1807--A1832.

\bibitem{Bachinietal_arXiv_2022}
E.~Bachini, P.~Brandner, T.~Jankuhn, M.~Nestler, S.~Praetorius, A.~Reusken,
  A.~Voigt, Diffusion of tangential tensor fields: numerical issues and
  influence of geometric properties, arXiv  2205.12581.

\bibitem{jankuhn2018}
T.~Jankuhn, M.~A. Olshanskii, A.~Reusken, {Incompressible fluid problems on
  embedded surfaces: Modeling and variational formulations}, Interf. Free
  Bound. {20} ({2018}) {353--377}.

\bibitem{barrett2008parametric}
J.~W. Barrett, H.~Garcke, R.~N{\"u}rnberg, On the parametric finite element
  approximation of evolving hypersurfaces in {$I\!\!R^3$}, J. Comput. Phys. 227
  (2008) 4281--4307.

\bibitem{praetorius2020dunecurvedgrid}
S.~Praetorius, F.~Stenger, {DUNE-CurvedGrid--A DUNE module for surface
  parametrization}, Arch. Num. Software 22 (2020) 1--22.

\bibitem{Dziuk_NM_2008}
G.~Dziuk, {Computational parametric Willmore flow}, Numer. Math. 111 (2008)
  55--80.

\bibitem{Reutheretal_MMS_2015}
S.~Reuther, A.~Voigt, The interplay of curvature and vortices in flow on curved
  surfaces, Multiscale Model. Sim. 13 (2015) 632--643.

\bibitem{Kobaetal_QAM_2017}
H.~Koba, C.~Liu, Y.~Giga, {Energetic variational approaches for incompressible
  fluid systems on an evolving surface}, Quart. Appl. Math. 75 (2017) 359--389.

\bibitem{Miura_QAM_2018}
T.-H. Miura, On singular limit equations for incompressible fluids in moving
  thin domains, Quart. Appl. Math. 76 (2018) 215--251.

\bibitem{Nitschkeetal_PRF_2019}
I.~Nitschke, S.~Reuther, A.~Voigt, Hydrodynamic interactions in polar liquid
  crystals on evolving surfaces, Phys. Rev. Fluids 4 (2019) 044002.

\bibitem{Brandneretal_arXiv_2021}
P.~Brandner, A.~Reusken, P.~Schwering, On derivations of evolving surface
  navier-stokes equations, arXiv:2110.14262 (2021).

\bibitem{Reutheretal_MMS_2018}
S.~Reuther, A.~Voigt, Erratum: "{T}he interplay of curvature and vortices in
  flow on curved surfaces", Multiscale Model. Sim. 16 (2018) 1448--1453.

\bibitem{Kobaetal_QAM_2018}
H.~Koba, C.~Liu, Y.~Giga, {Errata to “Energetic variational approaches for
  incompressible fluid systems on an evolving surface"}, Quart. Appl. Math. 76
  (2018) 174--152.

\bibitem{Yavarietal_JNS_2016}
A.~Yavari, A.~Ozakin, S.~Sadik, Nonlinear elasticity in a deforming ambient
  space, J. Nonli. Sci. 26 (2016) 1651--1692.

\bibitem{Helfrich+1973+693+703}
W.~Helfrich, Elastic properties of lipid bilayers: Theory and possible
  experiments, Zeitschrift für Naturforschung C 28 (1973) 693--703.

\bibitem{Dziuk_NM_1990}
G.~Dziuk, {An algorithm for evolutionary surfaces}, Numer. Math. 58 (1990)
  603--611.

\bibitem{Baenschetal_JCP_2005}
E.~Bänsch, P.~Morin, R.~H. Nochetto, A finite element method for surface
  diffusion: the parametric case, J. Comput. Phys. 203 (2005) 321--343.

\bibitem{Hausseretal_JSC_2007}
F.~Hau{\ss}er, A.~Voigt, A discrete scheme for parametric anisotropic surface
  diffusion, J. Sci. Comput. 30 (2007) 223--235.

\bibitem{demlow2006localized}
A.~Demlow, Localized pointwise a posteriori error estimates for gradients of
  piecewise linear finite element approximations to second-order quasilinear
  elliptic problems, SIAM J. Numer. Anal. 44 (2006) 494--514.

\bibitem{bastian2021dune}
P.~Bastian, M.~Blatt, A.~Dedner, N.-A. Dreier, C.~Engwer, R.~Fritze,
  C.~Gr{\"a}ser, C.~Gr{\"u}ninger, D.~Kempf, R.~Kl{\"o}fkorn, et~al., The
  {DUNE} framework: Basic concepts and recent developments, Comput. Math. Appl.
  81 (2021) 75--112.

\bibitem{dunealugrid:16}
M.~Alk{\"a}mper, A.~Dedner, R.~Kl{\"o}fkorn, M.~Nolte, {The DUNE-ALUGrid
  module.}, Arch. Numer. Software 4 (2016) 1--28.

\bibitem{Kovacsetal_NM_2019}
B.~Kovacs, B.~Li, C.~Lubich, {A convergent evolving finite element algorithm
  for mean curvature flow of closed surfaces}, Numer. Math. 143 (2019)
  797--853.

\bibitem{dziuk2013finite}
G.~Dziuk, C.~M. Elliott, Finite element methods for surface {PDEs}, Acta
  Numerica 22 (2013) 289--396.

\bibitem{Bonitoetal_JCP_2010}
R.~H. Bonito, A.~Nochetto, M.~S. Pauletti, Parametric fem for geometric
  biomembranes, J. Comput. Phys. 229 (2010) 3171--3188.

\bibitem{Nitschkeetal_PAMM_2020}
I.~Nitschke, S.~Reuther, A.~Voigt, Vorticity-stream function approaches are
  inappropriate to solve the surface navier-stokes equation on a torus, Proc.
  Appl. Math. Mech. 20 (2020) e202000006.

\bibitem{vey2007amdis}
S.~Vey, A.~Voigt, {AMDiS}: adaptive multidimensional simulations, Comput. Vis.
  Sci. 10 (2007) 57--67.

\bibitem{witkowski2015software}
T.~Witkowski, S.~Ling, S.~Praetorius, A.~Voigt, Software concepts and numerical
  algorithms for a scalable adaptive parallel finite element method, Adv.
  Comput. Math. 41 (2015) 1145--1177.

\bibitem{Seifert_AP_1997}
U.~Seifert, Configurations of fluid membranes and vesicles, Adv. Phys. 46
  (1997) 13--137.

\bibitem{Nitschkeetal_TFI_2017}
I.~Nitschke, S.~Reuther, A.~Voigt, Discrete exterior calculus {(DEC)} for the
  surface {Navier-Stokes} equation, in: D.~Bothe, A.~Reusken (Eds.), Transport
  Processes at Fluidic Interfaces, 2017, pp. 177--197.

\bibitem{Pruessetal_JEE_2021}
J.~Pr\"uss, G.~Simonett, M.~Wilke, On the {Navier–Stokes} equations on
  surfaces, J. Evol. Eq. 21 (2021) 3153--3179.

\bibitem{hughes1996application}
S.~W. Hughes, T.~J. D'Arcy, D.~J. Maxwell, J.~E. Saunders, C.~F. Ruff, W.~S.~C.
  Chiu, R.~J. Sheppard, Application of a new discret form of {Gauss'} theorem
  for measuring volume, Phys. Med. Bio. 41 (1996) 1809.

\bibitem{bao2021structure}
W.~Bao, Q.~Zhao, A structure-preserving parametric finite element method for
  surface diffusion, SIAM J. Numer. Anal. 59 (2021) 2775--2799.

\bibitem{Kovacs_IAMJNA_2018}
B.~Kovacs, High-order evolving surface finite element method for parabolic
  problems on evolving surfaces, IMA J. Numer. Anal. 38 (2018) 430--459.

\bibitem{Elliottetal_NM_2022}
C.~M. Elliott, H.~Garcke, B.~Kovacs, Numerical analysis for the interaction of
  mean curvature flow and diffusion on closed surfaces, Numer. Math. 151 (2022)
  873--925.

\bibitem{DLM_IMAJNA_2012}
G.~Dziuk, C.~Lubich, D.~Mansour, {Runge–Kutta} time discretization of
  parabolic differential equations on evolving surfaces, IMA J. Numer. Anal. 32
  (2012) 394--416.

\bibitem{LMV_IMAJNA_2013}
C.~Lubich, D.~Mansour, C.~Venkrataraman, Backward difference time
  discretization of parabolic differential equations on evolving surfaces, IMA
  J. Numer. Anal. 33 (2013) 1365--1385.

\bibitem{BR_IFB_2012}
N.~Balzani, M.~Rumpf, A nested variational time discretization for parametric
  willmore flow, Interf. Free Bound. 14 (2012) 413--454.

\bibitem{Boisetal_PRL_2011}
J.~S. Bois, F.~J\"ulicher, S.~W. Grill, Pattern formation in active fluids,
  Phys. Rev. Lett. 106 (2011) 028103.

\bibitem{Mickelinetal_PRL_2019}
O.~Mickelin, J.~S\l{}omka, K.~J. Burns, D.~Lecoanet, G.~M. Vasil, L.~M. Faria,
  J.~Dunkel, Anomalous chained turbulence in actively driven flows on spheres,
  Phys. Rev. Lett. 120 (2018) 164503.

\bibitem{Apaza_SM_2018}
L.~Apaza, M.~Sandoval, {Active matter on Riemannian manifolds}, {Soft Matter}
  {14} ({2018}) {9928--9936}.

\bibitem{Pearceetal_PRL_2019}
D.~Pearce, P.~Ellis, A.~Fernandez-Nieves, L.~Giomi, {Geometrical control of
  active turbulence in curved topographies}, {Phys. Rev. Lett.} {122} ({2019})
  {168002}.

\bibitem{Nitschkeetal_PRSA_2020}
I.~Nitschke, S.~Reuther, A.~Voigt, Liquid crystals on deformable surfaces,
  Proc. Roy. Soc. A 476 (2020) 20200313.

\bibitem{Al-Izzietal_SCDB_2021}
S.~C. Al-Izzi, R.~G. Morris, Active flows and deformable surfaces in
  development, Seminars in Cell \& Developmental Biology 120 (2021) 44--52.

\bibitem{Hoffmannetal_SA_2022}
L.~A. Hoffmann, L.~N. Carenza, J.~Eckert, L.~Giomi, {Theory of defect-mediated
  morphogenesis}, {Sci. Adv.} {8} ({2022}) {eabk2712}.

\end{thebibliography}

\end{document}